\documentclass[10pt]{amsart}
\usepackage{amsfonts}
\usepackage{graphicx}
\usepackage{amscd}
\usepackage{amsmath,amsfonts,amssymb,amsthm,mathrsfs,xypic}
\xyoption{curve}
\usepackage{fancybox}
\newcommand{\m}{\Lambda}
\newcommand{\Hom}{\operatorname{Hom}}
\newcommand{\End}{\operatorname{End}}
\newcommand{\Ext}{\operatorname{Ext}}
\newcommand{\Ker}{\operatorname{Ker}}
\newcommand{\cok}{\operatorname{Coker}}
\newcommand{\Ima}{\operatorname{Im}}

\newcommand{\s}{\hfill \blacksquare}
\newtheorem{thm}{Theorem}[section]

\newtheorem{cor}[thm]{Corollary}

\newtheorem{lem}[thm]{Lemma}

\newtheorem{prop}[thm]{Proposition}

\newtheorem{defn}[thm]{Definition}
\begin{document}

\title [Monomorphism operator and perpendicular operator] {Monomorphism operator and perpendicular operator}
\author [Keyan Song, Pu Zhang] {Keyan Song$^\dag$, Pu Zhang$^*$}
\thanks{$^*$ The corresponding author \ \ \ \ \ \ pzhang$\symbol{64}$sjtu.edu.cn}
\thanks{Supported by the NSF of China (11271251), and Doctoral Fund of Ministry of Education (20120073110058).}
\thanks {$^\dag$ skyysh$\symbol{64}$outlook.com}

\maketitle

\begin{center}
Department of Mathematics, \ \ Shanghai Jiao Tong University\\
Shanghai 200240, PR China
\end{center}

\begin{abstract} \ For a quiver $Q$,
a $k$-algebra $A$, and a full subcategory $\mathcal X$ of $A$-mod,
the monomorphism category ${\rm Mon}(Q, \mathcal X)$ is introduced.
The main result says that if $T$ is an $A$-module such that there is
an exact sequence $0\rightarrow T_m\rightarrow\cdots\rightarrow
T_0\rightarrow D(A_A)\rightarrow 0$ with each $T_i\in {\rm add}
(T)$, then ${\rm Mon}(Q, \ ^\perp T) = \ ^\perp (kQ\otimes_k T)$;
and if $T$ is cotilting, then $kQ\otimes_k T$ is a unique cotilting
$\m$-module, up to multiplicities of indecomposable direct summands,
such that ${\rm Mon}(Q, \ ^\perp T)= \ ^\perp (kQ \otimes_k T)$.

As applications, the category of the Gorenstein-projective
$(kQ\otimes_kA)$-modules is characterized as ${\rm Mon}(Q,
\mathcal{GP}(A))$ if $A$ is Gorenstein; the contravariantly
finiteness of ${\rm Mon}(Q, \mathcal X)$ can be described; and a
sufficient and necessary condition for ${\rm Mon}(Q, A)$ being of
finite type is given.

\vskip15pt


\vskip5pt

\end{abstract}

\section {\bf Introduction}

\subsection{} With a quiver $Q$ and a $k$-algebra $A$, one can associate the monomorphism category ${\rm
Mon}(Q, A)$ ([LZ]). If $Q = \bullet \rightarrow \bullet$ it is
called {\it the submodule category} and denoted by $\mathcal S(A)$.
If $Q = n \bullet \rightarrow \cdots \rightarrow
 \bullet 1$ it is called {\it the
filtered chain category} in D. Simson [S]; and it is denoted by
$\mathcal S_n(A)$ in [Z].

\vskip10pt

G. Birkhoff [B] initiates the study of $\mathcal S(\Bbb Z/\langle
p^t\rangle)$. C. M. Ringel and M. Schmidmeier ([RS1] - [RS3]) have
extensively studied $\mathcal S(A)$. In particular, the
Auslander-Reiten theory of $\mathcal S(A)$ is explicitly given
([RS2]). Since then the monomorphism category receives more
attention. In [Z] relations among  $\mathcal S_n(A)$ and the
Gorenstein-projective modules and cotilting theory are given. D.
Kussin, H. Lenzing, and H. Meltzer [KLM1] establish a surprising
link between the stable submodule category and the singularity
theory via weighted projective lines (see also [KLM2]). In [XZZ] the
Auslander-Reiten theory of $\mathcal S(A)$ is extended to $\mathcal
S_n(A)$. For more related works we refer to [A], [RW], [SW], [Mo],
[C1], [C2], and [RZ].

\vskip10pt

\subsection{} Let $\mathcal X$ be a full subcategory of $A$-mod. We also define the monomorphism category
${\rm Mon}(Q, \mathcal X)$. For an $A$-module $T$, let $^\perp T$ be
the full subcategory of $A$-mod consisting of those modules $X$ with
$\Ext_A^i(X, T)=0, \ \forall \ i\ge 1$. The main result of this
paper gives a reciprocity  of the monomorphism operator ${\rm
Mon}(Q, -)$ and the left perpendicular operator $^\bot$.  Namely, if
$T$ is an $A$-module such that there is an exact sequence
$0\rightarrow T_m\rightarrow\cdots\rightarrow T_0\rightarrow
D(A_A)\rightarrow 0$ with each $T_i\in {\rm add} (T)$, then ${\rm
Mon}(Q, \ ^\perp T) = \ ^\perp (kQ\otimes_k T)$ \ (Theorem
\ref{mainthm1}); and if $T$ is a cotilting $A$-module, then
$kQ\otimes_k T$ is a unique cotilting $\m$-module, up to
multiplicities of indecomposable direct summands, such that ${\rm
Mon}(Q, \ ^\perp T)= \ ^\perp (kQ \otimes_k T)$ \ (Theorem
\ref{mainthm2}).

\vskip10pt

Theorems \ref{mainthm1} and \ref{mainthm2} generalize [Z, Theorem
3.1${\rm (i)}$ and ${\rm (ii)}$] for $Q = \bullet \rightarrow \cdots
\rightarrow
 \bullet$. However, the arguments in [Z] can not be generalized to the general
case (cf. 3.1 and 4.1 below). Here we adopt new treatments, in
particular by using an adjoint pair $(\cok_i, S(i)\otimes -)$ and
Lemma \ref{self-orthogonal}.

\vskip10pt

\subsection{}  Our main results have some applications, which generalize the corresponding results in [Z].

\vskip10pt

The category $\mathcal {GP}(A)$ of the Gorenstein-projective
$A$-modules is Frobenius (cf. [AB], [AR], [EJ]), and hence the
corresponding stable category is triangulated ([H]). If $A$ is
Gorenstein (i.e., ${\rm inj.dim} _AA< \infty$ and ${\rm inj.dim} A_A
< \infty$), then $\mathcal {GP}(A) = \ ^\perp A$ ([EJ, Corollary
11.5.3]). Taking $T = \ _AA$ in Theorem \ref{mainthm1} we have
$\mathcal{GP}(\m) = {\rm Mon}(Q, \mathcal {GP}(A))$ if $A$ is
Gorenstein.

\vskip10pt

M. Auslander and I. Reiten [AR, Theorem 5.5(a)] have established a
deep relation between resolving contravariantly finite subcategories
and cotilting theory, by asserting that $\mathcal X$ is resolving
and contravariantly finite with $\widehat{\mathcal X} = A$-mod if
and only if $\mathcal X = \ ^\perp T$ for some cotilting $A$-module
$T$, where $\widehat{\mathcal X}$ is the full subcategory of $A$-mod
consisting of those modules $X$, such that there is an exact
sequence $0\rightarrow X_{m} \rightarrow \cdots \rightarrow
X_{0}\rightarrow X\rightarrow 0$ with each $X_{i}\in \mathcal{X}.$
It is natural to ask when is ${\rm Mon}(Q, \mathcal X)$
contravariantly finite in $\m$-mod? As an application of Theorem
\ref{mainthm2} and [AR, Theorem 5.5(a)], we see that
$\rm{Mon}(Q,\mathcal X)$ is resolving and contravariantly finite
with $\widehat{{\rm Mon}(Q, \mathcal X)}= \m$-mod if and only if
$\mathcal X$ is resolving and contravariantly finite with
$\widehat{\mathcal X} = A$-mod (Theorem \ref{app1}).

\vskip10pt

It is well-known that the representation type of ${\rm Mon}(Q, A)$
is different from the ones of $A$ and of $\m = kQ\otimes_k A$. For
example, $k[x]/\langle x^t\rangle$ is of finite type, while
$k(\bullet \rightarrow \bullet)\otimes_kk[x]/\langle x^t\rangle$ is
of finite type if and only if $t\le 3$, and $\mathcal
S_2(k[x]/\langle x^t\rangle)$ is of finite type if and only if $t\le
5$. If $t > 6$ then $\mathcal S_2(k[x]/\langle x^t\rangle)$ is of
``wild" type, while $\mathcal S_2(k[x]/\langle x^6\rangle)$ is of
``tame" type ([S], Theorems 5.2 and 5.5). A complete classification
of indecomposable objects of $\mathcal S_2(k[x]/\langle x^6\rangle)$
is exhibited in [RS3]. Inspired by Auslander's classical result: $A$
is representation-finite if and only if there is an
$A$-generator-cogenerator $M$ such that ${\rm gl.dim} \End_A(M)$
$\le 2$ ([Au], Chapter III), by using Theorem \ref{mainthm2} we
prove that ${\rm Mon}(Q, A)$ is of finite type if and only if there
is a generator and relative cogenerator $M$ of ${\rm Mon}(Q, A)$
such that ${\rm gl.dim} \End_{\m}(M)\le 2$ (Theorem \ref{app2}).

\vskip10pt

\section {\bf Preliminaries on monomorphism categories}

In this section we fix notations, and give necessary definitions and
facts.

\vskip10pt

\subsection{} Throughout this paper, $k$ is a field, $Q$ is a finite acyclic quiver (i.e., a finite quiver without oriented cycles), and
$A$ is a finite-dimensional $k$-algebra. Denote by $kQ$ the path
algebra of $Q$ over $k$. Put $\m = kQ\otimes_kA$, and $D = \Hom_k(-,
k).$ Let $P(i)$ (resp. $I(i)$) be the indecomposable projective
(resp. injective) $kQ$-module, and $S(i)$ the simple $kQ$-module, at
$i\in Q_0$. By $A$-mod we denote the category of finite-dimensional
left $A$-modules. For an $A$-module $T$, let ${\rm add}(T)$ be the
full the subcategory of $A$-mod consisting of all the direct sums of
indecomposable direct summands of $T$.

\vskip10pt

\subsection{} Given a finite acyclic quiver $Q = (Q_0,
Q_1, s, e)$ with $Q_0$ the set of vertices and $Q_1$ the set of
arrows, we write the conjunction of a path $p$ of $Q$ from right to
left, and let $s(p)$ and $e(p)$ be respectively the starting and the
ending point of $p$. The notion of representations of $Q$ over $k$
can be extended as follows. By definition ([LZ]), {\it a
representation $X$ of $Q$ over $A$} is a datum $X = (X_i, \
X_{\alpha}, \ i\in Q_0, \ \alpha \in Q_1)$, or simply $X = (X_i,
X_\alpha)$, where each $X_i$ is an $A$-module, and each $X_{\alpha}:
X_{s(\alpha)} \rightarrow X_{e(\alpha)}$ is an $A$-map. It is {\it a
finite-dimensional representation} if so is each $X_i$. We call
$X_i$ {\it the $i$-th branch} of $X$. A morphism $f$ from $X$ to $Y$
is a datum $(f_i, \ i\in Q_0)$, where $f_i: X_i \rightarrow Y_i$ is
an $A$-map for $i\in Q_0$, such that for each arrow $\alpha:
j\rightarrow i$ the following diagram

\[\xymatrix {X_j\ar[r]^{f_j}\ar[d]_{X_\alpha} & Y_j\ar[d]^{Y_\alpha} \\
X_i\ar[r]^{f_i} & Y_i}\eqno(2.1)\] commutes. Denote by ${\rm Rep}(Q,
A)$ the category of finite-dimensional representations of $Q$ over
$A$. Note that a sequence of morphisms $0\longrightarrow X\stackrel
{f}\longrightarrow Y\stackrel {g}\longrightarrow Z \longrightarrow
0$ in ${\rm Rep}(Q, A)$ is exact if and only if $0\longrightarrow
X_i\stackrel {f_i}\longrightarrow Y_i\stackrel {g_i}\longrightarrow
Z_i \longrightarrow 0$ is exact in $A$-mod for each $i\in Q_0$.

\vskip10pt

\begin{lem} \label{rep} {\rm ([LZ, Lemma 2.1])} \ We have an equivalence $\Lambda\mbox{-}{\rm
mod}\cong {\rm Rep}(Q, A)$ of categories.\end{lem}

In the following we will identify a $\Lambda$-module with a
representation of $Q$ over $A$. If $T\in A$-mod and $M\in kQ$-mod
with $M = (M_i, i\in Q_0, M_\alpha, \alpha\in Q_1)\in {\rm Rep}(Q,
k)$, then $M\otimes_k T\in \m$-mod with $M\otimes_k T =
(M_i\otimes_kT = T^{{\rm dim}_kM_i}, \ i\in Q_0, \
M_\alpha\otimes_k{\rm Id}_T, \ \alpha\in Q_1)\in {\rm Rep}(Q, A).$

\vskip10pt

 \subsection{} Here is the central notion of this
paper.
\begin{defn} \label{maindef} \
${\rm (i) \ \ ([LZ])}$ \ A representation $X = (X_i, X_{\alpha},
i\in Q_0, \alpha\in Q_1)\in {\rm Rep}(Q, A)$ is a monic
representation of $Q$ over $A$, or a monic $\Lambda$-module, if
$\delta_i(X)$ is an injective $A$-map for each $i\in Q_0$, where
$$\delta_i(X)=(X_{\alpha})_{\alpha\in Q_1, \ e(\alpha) = i}: \ \bigoplus\limits_{\begin
{smallmatrix} \alpha\in Q_1\\ e(\alpha) = i \end{smallmatrix}}
X_{s(\alpha)} \longrightarrow X_i.$$

Denote by ${\rm Mon}(Q, A)$ the full subcategory of ${\rm Rep}(Q,
A)$ consisting of all the monic representations of $Q$ over $A$,
which is called  the monomorphism category of $A$ over $Q$.

\vskip10pt

${\rm (ii)}$ \ Let $\mathcal X$ be a full subcategory of $A$-mod.
Denote by ${\rm Mon}(Q, \mathcal X)$ the full subcategory of ${\rm
Mon}(Q, A)$ consisting of all the monic representations $X = (X_i,
X_{\alpha})$, such that $X_i\in \mathcal X$ and $\cok\delta_i(X) \in
\mathcal X$ for all $i\in Q_0$. We call ${\rm Mon}(Q, \mathcal X)$
the monomorphism category of $\mathcal X$ over $Q$.
\end{defn}

\vskip10pt

If $\mathcal X = A$-mod then ${\rm Mon}(Q, \mathcal X) = {\rm
Mon}(Q, A)$. For $M\in kQ$-mod and $T\in A$-mod, it is clear that if
$M\in {\rm Mon}(Q, k)$ then $M\otimes_kT\in {\rm Mon}(Q, A)$. In
particular, $P(i)\otimes_k T\in {\rm Mon}(Q, A)$ for each $i\in
Q_0$.

\vskip10pt

Note that $D(\m_\m) \cong D(kQ_{kQ})\otimes_k D(A_A)$ as left
$\m$-modules. We need the following fact.

\vskip10pt

\begin{lem}\label{inj} {\rm ([LZ, Proposition 2.4])} \ Let ${\rm Ind}\mathcal{P}(A)$ (resp. ${\rm Ind} \mathcal{I}(A))$ denote the set of pairwise
non-isomorphic indecomposable projective (resp. injective)
$A$-modules. Then

$$ {\rm Ind} \mathcal{P}(\m) = \{P(i)\otimes_k
P \ | \ i\in Q_0, \  P\in {\rm Ind} \mathcal{P}(A)\}\subseteq {\rm
Mon}(Q, A),$$ and
$${\rm Ind} \mathcal{I}(\m) = \{I(i)\otimes_k
I \ | \ i\in Q_0, \  I\in {\rm Ind} \mathcal{I}(A)\}.$$

\vskip5pt

\noindent In particular, for $M\in kQ$-mod we have ${\rm proj.dim}
(M\otimes_k A) \le 1$, and \ ${\rm inj.dim} (M\otimes_k D(A_A)) \le
1$.
\end{lem}

\vskip10pt

\subsection{} Given $X = (X_j, X_\alpha)\in \m$-mod, for each $i\in Q_0$ we have
functors $F_i$ and $F^+_i$ from $\m\mbox{-}{\rm mod}$ to
$A\mbox{-}{\rm mod}$, respectively induced by $F_i(X) = X_i$ and
$F^+_i(X): = \bigoplus\limits_{\begin {smallmatrix} \alpha\in Q_1\\
e(\alpha) = i \end{smallmatrix}} X_{s(\alpha)}$ (if $i$ is a source
then $F^+_i(X): = 0$). \vskip10pt

We write $\cok\delta_i(X)$ (cf. Definition \ref{maindef} ${\rm
(i)}$) as $\cok_i(X)$. Then we have a functor $\cok_i:
\m\mbox{-}{\rm mod}\longrightarrow A\mbox{-}{\rm mod},$
explicitly given by $\cok_i(X): = X_i/\sum\limits_{\begin {smallmatrix} \alpha\in Q_1\\
e(\alpha) = i \end{smallmatrix}} \Ima X_{\alpha}$ (if $i$ is a
source then $\cok_i(X): = X_i$). So we have an exact sequence of
functors $F^{+}_i \stackrel {\delta_i}\longrightarrow F_i \stackrel
{\pi_i}\longrightarrow \cok_i \longrightarrow 0,$ i.e., we have the
exact sequence of $A$-modules $$\bigoplus\limits_{\begin
{smallmatrix} \alpha\in Q_1\\ e(\alpha) = i \end{smallmatrix}}
X_{s(\alpha)}\stackrel {\delta_i(X)}\longrightarrow  X_i \stackrel
{\pi_i(X)}  \longrightarrow \cok_i(X) \longrightarrow 0$$ for each
$X\in \m$-mod, where $\pi_i(X)$ is the canonical map. It is clear
that $F^{+}_i$ and $F_i$ are exact, and $\cok_i$ is right exact (by
Snake Lemma). For $i, j\in Q_0$ and $T\in A$-mod, we have
$$\cok_i(P(j)\otimes_k T)=\begin{cases} T, & \mbox{if} \ j=i;
\\ 0, & \mbox{if} \ j\ne i.
\end{cases} \eqno (2.2) $$

\vskip10pt

\begin{lem} \label{exact} \ For each $i\in Q_0$, the restriction of functor $\cok_i$ to ${\rm Mon}(Q,A)$ is exact.
\end{lem}

\noindent{\bf Proof.} Let
$0\rightarrow(X_i,X_\alpha)\rightarrow(Y_i,Y_\alpha)\rightarrow(Z_i,Z_\alpha)\rightarrow
0$ be an exact sequence in ${\rm Mon}(Q,A)$. Then we have the
following commutative diagram with exact rows
\[
\xymatrix{0 \ar[r] & \bigoplus X_{s(\alpha)} \ar[r]
\ar[d]^{\delta_i(X)} & \bigoplus Y_{s(\alpha)} \ar[r]
\ar[d]^{\delta_i(Y)} & \bigoplus Z_{s(\alpha)} \ar[d]^{\delta_i(Z)} \ar[r] & 0 \\
0\ar[r] & X_i \ar[r] & Y_i \ar[r] &Z_i \ar[r] & 0.}
\]
Then the assertion follows from Snake Lemma since $\delta_i(Z)$ is
injective. $\s$

\vskip10pt

Recall from [AR] that $\mathcal X$ is {\it resolving} if $\mathcal
X$ contains all the projective $A$-modules, $\mathcal X$ is closed
under taking extensions, kernels of epimorphisms, and direct
summands. Dually one has {\it a coresolving subcategory}.

\vskip10pt

\begin{lem} \label {exact2} \  Let $\mathcal X$ be a full subcategory of $A$-mod. Then

\vskip5pt

${\rm (i)}$ \ ${\rm Mon}(Q, \mathcal X)$ is closed under taking
extensions (resp. kernels of epimorphisms, direct summands) if and
only if $\mathcal X$ is closed under taking extensions (resp.
kernels of epimorphisms, direct summands).

\vskip5pt

${\rm (ii)}$ \ ${\rm Mon}(Q, \mathcal X)$ is resolving if and only
if  $\mathcal X$ is resolving. In particular, ${\rm Mon}(Q, A)$ is
resolving.
\end{lem}
\noindent {\bf Proof.} ${\rm (i)}$ can be similarly proved
 as  Lemma \ref{exact}. For ${\rm (ii)}$, \ by Lemma \ref{inj}
the branches of projective $\m$-modules are projective $A$-modules.
From this and ${\rm (i)}$ the assertion follows. $\s$

\vskip10pt

 \section{\bf Reciprocity}

\subsection{} This section is to prove the following reciprocity of the
monomorphism operator and the left perpendicular operator.

\vskip10pt

\begin{thm} \label{mainthm1} \ Let $T$ be an
$A$-module such that there is an exact sequence $0\rightarrow
T_m\rightarrow  \cdots \rightarrow T_0\rightarrow D(A_A)\rightarrow
0$ with each $T_j\in {\rm add} (T)$, then ${\rm Mon}(Q, \ ^\perp T)
= \ ^\perp (kQ\otimes_k T).$
\end{thm}

\vskip10pt

For $Q = \bullet \rightarrow \cdots \rightarrow
 \bullet$ this result has been obtained in [Z, Theorem 3.1${\rm (i)}$].
Since some adjoint pairs in [Z, Lemma 1.2] are not available here,
the arguments in [Z] can not be generalized to the general case.
Here we adopt the following adjoint pair $(\cok_i, S(i)\otimes -)$.

\vskip10pt

\subsection{}The following observation will be used throughout this section.
\vskip10pt

\begin{lem} \label{adj} \ Let $X=(X_i, X_\alpha) \in \Lambda$-mod and  $T \in A$-mod. Then for each $i\in Q_0$
we have an isomorphism of abelian groups  which is natural in both
positions
$${\rm Hom}_A(\cok_i(X), T)\cong{\rm Hom}_\Lambda(X, S(i)\otimes_k T).$$
\end{lem}
\noindent{\bf Proof.} \ If we write $S(i)\otimes_k T\in \m$-mod as
$(Y_j, Y_\alpha)$, then $Y_j = 0$ for $j\ne i$ and $Y_i = T$.
Consider  the homomorphism $\Psi: {\rm Hom}_A(\cok_i(X),
T)\rightarrow {\rm Hom}_\Lambda(X, S(i)\otimes_k T)$ given by
$$f\mapsto \Psi(f) = (g_j, j\in Q_0): \ X\rightarrow S(i)\otimes_kT,
\ \forall f\in {\rm Hom}_A(\cok_i(X), T),$$ where $g_j = 0$ for
$j\ne i$, and $g_i = f \ \pi_i(X): X_i \rightarrow T$ with the
canonical map $\pi_i(X): X_i\rightarrow \cok_i(X)$. By (2.1) it is
clear that $\Psi(f)\in {\rm Hom}_\Lambda(X, S(i)\otimes_k T)$ and
$\Psi$ is surjective. It is injective since $\pi_i(X)$ is
surjective. $\s$

\vskip10pt

\subsection{} \ We need the following fact.

\vskip10pt

\begin{lem} \label{perp1}  Let $T$ be an $A$-module. For each $i\in Q_0$ we have \ $^\perp (kQ\otimes_k T) = \ ^\perp (\bigoplus\limits_{i\in
Q_0} (S(i)\otimes_k T)).$
\end{lem}
\noindent{\bf Proof.} \ Put $S = \bigoplus\limits_{i\in Q_0} S(i)$,
and $J$ to be the Jacobson radical of $kQ$ with $J^l = 0$. Let $X\in
\ ^\perp (kQ\otimes T)$. By the exact sequence $0\rightarrow
J\otimes_k T\rightarrow kQ\otimes_k T\rightarrow S\otimes_k
T\rightarrow 0$ we get the exact sequence:
 $$\cdots\rightarrow  {\rm Ext}_\Lambda^j(X, kQ\otimes_k T)\rightarrow {\rm Ext}^j_\Lambda(X, S\otimes_k T)
 \rightarrow {\rm Ext}^{j+1}_\Lambda(X, J\otimes_k T)
 \rightarrow\cdots.$$
Since $kQ$ is hereditary, $J\in {\rm add} (kQ)$ and hence ${\rm
Ext}^j_\Lambda (X, J\otimes_k T)= 0, \ \forall \ j\ge 1$. Thus $X\in
\ ^\perp (S\otimes_k T)$.

Conversely, let  $X\in \ ^\perp (S\otimes_k T)$. From the exact
sequence $0\rightarrow J^{l-1}\otimes_kT \rightarrow
J^{l-2}\otimes_kT \rightarrow (J^{l-2}/J^{l-1})\otimes_kT\rightarrow
0$ and by $J^{l-1}\otimes_kT, \ J^{l-2}/J^{l-1}\otimes_kT\in {\rm
add}(S\otimes_k T)$ we see $X\in \ ^\perp (J^{l-2}\otimes_k T)$.
Continuing this process we finally see $X\in \ ^\perp
(J^{0}\otimes_k T) = \ ^\perp (kQ\otimes_k T)$. $\s$

 \vskip10pt

\begin{prop}\label{perp2}  \ We have ${\rm Mon}(Q,A)= \ ^\perp (kQ\otimes_k D(A_A))$.

\end{prop} \noindent{\bf Proof.}  By Lemma
\ref{perp1} it suffices to prove the following equality, for each
$i\in Q_0$:

$$^\perp (S(i)\otimes_k D(A_A)) = \{X=(X_j, X_\alpha)\in {\rm Rep}(Q, A) \ | \
\delta_i(X)\ \text{is injective}\}. $$

\vskip5pt

Let $P_X$ be the projective cover of $X$. Applying functor $F_i^+$
and $F_i$ to the exact sequence $0\rightarrow \Omega(X)\rightarrow
P_X\rightarrow X\rightarrow 0$ we get the following commutative
diagram with exact rows
\[\begin{CD}
0 @>>>F^{+}_i(\Omega(X))@>>>F^{+}_i(P_X)@>>>F^{+}_i(X)@>>>0\\
@. @VV\delta_i(\Omega(X))V @VV\delta_i(P_X)V @VV\delta_i(X)V \\
0 @>>>F_i(\Omega(X))@>>>F_i(P_X)@>>>F_i(X)@>>>0.\\
\end{CD}\]
By Snake Lemma we have the  exact sequence
$$0\rightarrow \Ker\delta_i(X)\rightarrow \cok_i(\Omega(X))\rightarrow \cok_i(P_X)\rightarrow \cok_i(X)\rightarrow 0.\eqno(*)$$

Assume that $\delta_i(X)$ is injective. Applying ${\rm Hom}_A(-,
D(A_A))$ to $(*)$ and by Lemma \ref{adj} we get the following exact
sequence (with Hom omitted)
$$0\rightarrow (X, S(i)\otimes_k D(A_A))\rightarrow (P_X, S(i)\otimes_k
D(A_A))\rightarrow (\Omega(X), S(i)\otimes_k D(A_A))\rightarrow
0.\eqno(**)$$ Applying ${\rm Hom}_\m(-, S(i)\otimes_k D(A_A))$ to
$0\rightarrow \Omega(X)\rightarrow P_X\rightarrow X\rightarrow 0$ we
get the exact sequence
$$0\rightarrow (X, S(i)\otimes D(A_A))\rightarrow (P_X, S(i)\otimes D(A))\rightarrow (\Omega(X), S(i)\otimes D(A))
\rightarrow {\rm Ext}^1_\m(X, S(i)\otimes D(A))\rightarrow 0.$$
Comparing it with $(**)$ we see $ {\rm Ext}^1_\m(X, S(i)\otimes_k
D(A_A))= 0$. By Lemma \ref {inj} we have ${\rm
inj.dim}(S(i)\otimes_k D(A_A))\leq 1$, so $X\in \ ^\perp
(S(i)\otimes_k D(A_A))$.

\vskip5pt

Conversely, assume $X\in \ ^\perp (S(i)\otimes_k D(A_A))$. Applying
${\rm Hom}_\m(-, S(i)\otimes_k D(A_A))$ to $0\rightarrow
\Omega(X)\rightarrow P_X\rightarrow X\rightarrow 0$ and using Lemma
\ref{adj}, we get the following exact sequence
$$0 \rightarrow (\cok_i(X), D(A_A)) \rightarrow (\cok_i(P_X), D(A_A)) \rightarrow (\cok_i(\Omega(X)), D(A_A))\rightarrow
0,$$ i.e., $0\rightarrow \cok_i(X)\rightarrow
\cok_i(P(X))\rightarrow \cok_i(\Omega(X))\rightarrow0$ is exact.
Comparing it with $(*)$ we see $\Ker\delta_i(X)=0$. $\s$

\vskip10pt

\subsection{} Replacing $D(A_A)$ in Proposition
\ref{perp2} by an arbitrary $A$-module $T$, we have

\vskip10pt

\begin{prop} \label{perp3} \ Let $T$ be an $A$-module. Then ${\rm Mon}(Q, \ ^\perp
T) = \ ^\perp (kQ\otimes_k T)\ \cap\ {\rm Mon}(Q, A)$.
\end{prop}

\noindent{\bf Proof.} We first prove that for each $i\in Q_0$ there
holds the following equality \begin{align*} \ \ \ \ \ ^\perp
(S(i)\otimes_k T) \ \cap \ {\rm Mon}(Q, A)  = \{ & X =(X_j,
X_\alpha)\in {\rm Rep}(Q, A) \ | \ \cok_i(X)\in \ ^\perp T, \\  & \
\ \delta_j (X)\ \text{is injective for all} \ j\in Q_0\}. \ \ \ \ \
\ \ \ \ \ \ \ \ \ \ \ \ \ \ \ \ \ \ \ \ \ \ \ \  (3.1)\end{align*}

Let $X\in {\rm Mon}(Q, A)$ with a projective resolution
$\cdots\rightarrow P^1\rightarrow P^0\rightarrow X\rightarrow 0$.
Since each $P^i$ is in ${\rm Mon}(Q, A)$ (cf. Lemma \ref{inj}) and
${\rm Mon}(Q, A)$ is closed under taking the kernels of epimorphisms
(cf. Lemma \ref{exact2}), it follows from Lemma \ref{exact} that we
have the exact sequence
$$\cdots\rightarrow \cok_i(P^1)\rightarrow \cok_i(P^0)\rightarrow
\cok_i(X)\rightarrow 0.$$ We claim it is a projective resolution of
$\cok_i(X)$. In fact, by $(2.2)$ we have
$$\cok_i(P(j)\otimes_k T)=\begin{cases} T, & \mbox{if} \ j=i;
\\ 0, & \mbox{if} \ j\ne i.
\end{cases} $$
So $\cok_i(kQ\otimes_k T)= T$ and $\cok_i(kQ\otimes_k A)= A $. Thus
$ \cok_i(P^j)$ is a projective $A$-module since $P^j\in {\rm add}
(kQ\otimes_k A).$

Applying ${\rm Hom}(-, S(i)\otimes_k T)$ to $\cdots\rightarrow
P^1\rightarrow P^0\rightarrow X\rightarrow 0$, by Lemma \ref{adj} we
have the following commutative diagram

\[\begin{CD}
0 @>>> (X, S(i)\otimes_k T) @>>>(P^0, S(i)\otimes_k T)@>>>(P^1, S(i)\otimes_k T)@>>>\cdots\\
@.  @ V   V \wr V @V V \wr V @ V V \wr V  \\
0 @>>>(\cok_i(X), T) @>>>(\cok_i(P^0), T) @>>>(\cok_i(P^1), T)@>>>\cdots\\
\end{CD}\]

\vskip5pt \noindent Note that  $X\in \ ^\perp (S(i)\otimes_k T)$ if
and only if the upper row is exact, if and only if the lower one is
exact, if and only if $\cok_i(X)\in  \ ^\perp T$. This proves
$(3.1)$.

\vskip10pt

Now, assume that $X\in {\rm Mon}(Q, \ ^\perp T)$. By definition and
$(3.1)$ we know $X\in \ ^\perp (S(i)\otimes_k T) \ \cap \ {\rm
Mon}(Q, A)$ for each $i\in Q_0$. By Lemma \ref{perp1} we know $X\in
\ ^\perp (kQ\otimes_k T)$ and hence  $X\in \ ^\perp (kQ\otimes_k T)
\ \cap \ {\rm Mon}(Q, A)$.

\vskip5pt

Conversely, assume that $X\in \ ^\perp (kQ\otimes_k T) \ \cap \ {\rm
Mon}(Q, A)$. By Lemma \ref{perp1} $X\in \ ^\perp (S(i)\otimes_k T) \
\cap \ {\rm Mon}(Q, A)$ for each $i\in Q_0$. To see $X\in {\rm
Mon}(Q, \ ^\perp T)$, by $(3.1)$ it remains to prove $X_i\in \
^\perp T$ for each $i\in Q_0$. For each $i\in Q_0$, set $l_i = 0$ if
$i$ is a source, and $l_i = {\rm max} \{ \ l(p) \ \mid p \ \mbox{is
a path with} \ e(p) = i\}$ if otherwise, where $l(p)$ is the length
of $p$. We prove $X_i\in \ ^\perp T$ by using induction on $l_i$. If
$l_i = 0$, then $i$ is a source and $X_i=\cok_i(X)\in \ ^\perp T$.
Let $l_i\ne 0$. Then we have the exact sequence $0 \longrightarrow
\bigoplus\limits_{\begin {smallmatrix} \alpha\in Q_1\\ e(\alpha) = i
\end{smallmatrix}} X_{s(\alpha)} \longrightarrow X_i \longrightarrow
\cok_i(X)\longrightarrow 0$ with $\cok_i(X)\in \ ^\perp T$.  Since
$l_{s(\alpha)} < l_i$ for $\alpha\in Q_1$ and $e(\alpha) = i$, by
induction $X_{s(\alpha)}\in \ ^\perp T$, and hence $X_i\in \ ^\perp
T$.  This completes the proof. $\s$

\vskip10pt

\subsection{Proof of Theorem \ref{mainthm1}.} \
By Proposition \ref{perp3} it suffices to prove $^\perp (kQ\otimes_k
T) \subseteq {\rm Mon}(Q, \ A)$. By Proposition \ref{perp2} it
suffices to prove $^\perp (kQ\otimes_k T) \subseteq  \ ^\perp
(kQ\otimes_k D(A_A))$. Let $X\in \ ^\perp (kQ\otimes_k T)$. By
assumption we have an exact sequence $0\longrightarrow kQ\otimes_k
T_m\longrightarrow \cdots \longrightarrow  kQ\otimes_k
T_0\longrightarrow kQ\otimes_k D(A_A)\longrightarrow 0$ with each
$kQ\otimes_k T_j \in {\rm add}(kQ\otimes_k T)$. From this we see the
assertion. $\s$

\vskip10pt

\subsection{} Let $\mathcal {GP}(A)$ denote the category  of the Gorenstein-projective
$A$-modules.  If $A$ is Gorenstein (i.e., ${\rm inj.dim} _AA<
\infty$ and ${\rm inj.dim} A_A < \infty$), then $\mathcal {GP}(A) =
\ ^\perp A$ ([EJ, Corollary 11.5.3]). Note that if $A$ is Gorenstein
then so is $\m$. Taking $T = \ _AA$ in Theorem \ref{mainthm1} we
have

\vskip10pt

\begin{cor} \label{gor} \ Let $A$ be a Gorenstein algebra. Then
$\mathcal{GP}(\m) = {\rm Mon}(Q, \mathcal {GP}(A))$.
\end{cor}

\vskip10pt

\section {\bf Monomorphism categories and cotilting theory}

\subsection{} The aim of this section is to prove the following

\vskip10pt

\begin{thm} \label {mainthm2} \ Let $T$ be a cotilting $A$-module.
Then $kQ\otimes_k T$ is a unique cotilting $\m$-module, up to
multiplicities of indecomposable direct summands, such that ${\rm
Mon}(Q, \ ^\perp T)= \ ^\perp (kQ \otimes_k T)$.
\end{thm}

\vskip10pt

For $Q = \bullet \rightarrow \cdots \rightarrow
 \bullet$ this result has been obtained in [Z, Theorem 3.1${\rm (ii)}$]. We
stress that the proof in [Z] can not be generalized to the general
case. Here we need to use Lemma \ref {self-orthogonal} below, rather
than a concrete construction in [Z, Lemma 3.7].

\vskip10pt

\subsection{} Recall that an $A$-module $T$ is an
$r$-cotilting module ([HR], [AR], [H], [Mi]) if the following
conditions are satisfied:

\vskip5pt

${\rm (i)}$ \ ${\rm inj.dim} T\leq r$;

\vskip5pt

${\rm (ii)}$ \  Ext$_A^i(T,T)=0$ for $i\geq 1$;

\vskip5pt

${\rm (iii)}$ \ there is an exact sequence $0\rightarrow
T_m\rightarrow  \cdots \rightarrow T_0\rightarrow D(A_A)\rightarrow
0$ with each $T_j\in {\rm add} (T)$.

\vskip10pt

For short, by ${\bf m}_i$ we denote the functor $P(i)\otimes_k- : \
A\mbox{-}{\rm mod} \rightarrow {\rm Mon}(Q,A)$, and by ${\bf m}$ we
denote the functor $kQ\otimes_k-: A\mbox{-}{\rm mod} \rightarrow
{\rm Mon}(Q,A)$. Then ${\bf m}(T) = \bigoplus\limits_{\begin
{smallmatrix} i\in Q_0
\end{smallmatrix}} {\bf m}_i(T) = kQ\otimes_k T, \ \forall \ T\in A\mbox{-}{\rm mod}.$

\vskip10pt

\begin{lem} \label{adj2} {\rm ([LZ, Lemma 2.3])} \ We have adjoint pair $({\bf m}_i, F_i)$ for each $i\in Q_0$, where functor
$F_i$ is defined in 2.4.
\end{lem}

\vskip10pt

We also need the following fact.

\vskip10pt

\begin{lem} \label{adj3} \  Let   $X=(X_j, X_\alpha) \in \m$-mod and $T \in A$-mod. Then we have an isomorphism of
abelian groups for each $i\in Q_0$, which is natural in both
positions
$${\rm Ext}_\m^s({\bf m}_i(T), X)\cong {\rm Ext}_A^s(T, X_i), \ \ \forall \ s\geq 0.$$
\end{lem}
\noindent{\bf Proof.} The proof is same as in [Z, Lemma 3.4] for $Q
= \bullet \rightarrow \cdots \rightarrow
 \bullet$. For completeness we
include a justification. Taking the $i$-th branch of an injective
resolution $0\rightarrow X\rightarrow I^0\rightarrow I^1\rightarrow
\cdots$ of $_\m X$, by Lemma \ref{inj} $0\rightarrow X_i\rightarrow
I_i^0\rightarrow I_i^1\rightarrow \cdots$ is an injective resolution
of $_AX_i$. On the other hand by Lemma \ref{adj2} we have the
following commutative diagram
\[\begin{CD}
0 @>>>{\rm Hom}_\m({\bf m}_i(T), X)@>>> {\rm Hom}_\m({\bf m}_i(T), I^0) @>>>{\rm Hom}_\m({\bf m}_i(T), I^1)@>>>\cdots\\
@.  @ V   V \wr V @V V \wr V @ V V \wr V  \\
0 @>>>{\rm Hom}_A(T, X_i)@>>> {\rm Hom}_A(T, I_i^0) @>>>{\rm Hom}_A(T, I_i^1)@>>>\cdots,\\
\end{CD}\]
from this we see the assertion. $\s$

\vskip10pt

\subsection {} Let $\mathcal X$ be a full subcategory of $A$-mod. Following
[AR] let $\widehat{\mathcal X}$ denote the full subcategory of
$A$-mod consisting of those $A$-modules $X$ such that there is an
exact sequence $0\rightarrow X_m\rightarrow X_{m-1}\rightarrow
\cdots\rightarrow X_0\rightarrow X\rightarrow 0$ with each $X_i\in
\mathcal X$. Recall that $\mathcal X$ is {\it self-orthogonal} if
$\Ext_A^s(M, N)=0, \ \ \forall \ M, \ N\in \mathcal X,  \ \forall \
s\geq 1.$ In this case $\widehat{\mathcal X} \subseteq \mathcal
X^\perp$, where $\mathcal X^\perp = \{X\in A\mbox{-}{\rm mod} \ | \
\Ext_A^i(M, X) = 0, \ \forall \ M\in \mathcal X, \ \forall \ i\ge
1\}$.

\vskip10pt

The following fact is of independent interest. It is a key step in
the proof of Theorem \ref{mainthm2}.

\vskip10pt

\begin{lem} \label {self-orthogonal} \  Let $\mathcal X$ be a self-orthogonal full
subcategory of $A$-mod. Then

\vskip5pt

${\rm (i)}$ \ $\widehat{\mathcal X}$ is closed under taking
cokernels of monomorphisms.

\vskip5pt

${\rm (ii)}$ \ $\widehat{\mathcal X}$ is closed under taking
extensions.

\vskip5pt

${\rm (iii)}$ \ If $\mathcal X$ is closed under taking kernels of
epimorphisms, then so is $\widehat{\mathcal X}$.

\end{lem}

\noindent{\bf Proof.} \  ${\rm (i)}$ \ Let $0 \rightarrow X
\xrightarrow {f} Y \rightarrow  Z \rightarrow 0$ be an exact
sequence with $X, Y\in \widehat{\mathcal X}$. By definition there
exist exact sequences $0\rightarrow X_n\rightarrow
X_{n-1}\rightarrow \cdots \rightarrow X_0\xrightarrow {c_0}
X\rightarrow 0$, and $0\rightarrow Y_n\rightarrow Y_{n-1}\rightarrow
\cdots \rightarrow Y_0\xrightarrow {d_0}Y\rightarrow 0$ with $X_i, \
Y_i\in \mathcal X \cup \{0\}, \ 0\le i\le n$. Since $\mathcal X$ is
self-orthogonal, $f: X\rightarrow Y$ induces a chain map $f^\bullet:
X^\bullet\rightarrow Y^\bullet$, where $X^\bullet$ is the complex
$0\rightarrow X_n\rightarrow X_{n-1}\rightarrow \cdots \rightarrow
X_0\rightarrow0$, and similarly for $Y^\bullet$. Consider the
following commutative diagram in the bounded derived category
${D}^b(A)$, with two rows being distinguished triangles
\[
\xymatrix{X^\bullet \ar[r]^{f^\bullet} \ar[d]^{c_0} &Y^\bullet
\ar[r] \ar[d]^{d_0}
&{\rm Con}(f^\bullet) \ar @{.>}[d] \ar[r]& X^\bullet [1] \ar[d] \\
X \ar[r] &Y \ar[r] & Z \ar[r] & X[1]}
\]

\vskip5pt \noindent (note that the lower row is also a distinguished
triangle since $0 \rightarrow X \xrightarrow {f} Y \rightarrow  Z
\rightarrow 0$ is exact), where ${\rm Con}(f^\bullet)$ is the
mapping cone $0\rightarrow X_n\rightarrow X_{n-1}\oplus
Y_n\rightarrow \cdots\rightarrow X_0\oplus Y_1\xrightarrow
{\partial} Y_0\rightarrow 0$. Since $c_0$ and $d_0$ are isomorphisms
in ${D}^b(A)$, we have $Z\cong {\rm Con}(f^\bullet)$ in ${D}^b(A)$.
It follows that the $i$-th cohomology group of ${\rm
Con}(f^\bullet)$ is isomorphic to the $i$-th cohomology group of the
stalk complex $Z$ for each $i\in \Bbb Z$. In particular ${\rm
Con}(f^\bullet)$ is exact except at the $0$-th position, and
$Y_0/\Ima\partial \cong Z$. Thus
$$0\rightarrow X_n\rightarrow X_{n-1}\oplus Y_n\rightarrow
\cdots\rightarrow X_0\oplus Y_1\stackrel{\partial}\longrightarrow
 Y_0 \rightarrow Z\rightarrow 0$$ is exact. This proves $Z\in \widehat{\mathcal X}$.

\vskip5pt

${\rm (iii)}$ can be similarly proved, and  ${\rm (ii)}$ can be
proved by a version of Horse-shoe Lemma. We omit the details. (Only
${\rm (i)}$ will be needed in the proof of Theorem \ref{mainthm2}.)
\ $\s$

\vskip10pt

\begin{lem} \label{cotilting} \ Let $T$ be an $r$-cotilting $A$-module. Then
$kQ\otimes_k T$ is an $(r+1)$-cotilting $\m$-module with $\End_\m
(kQ\otimes_k T)\cong (kQ\otimes_k \End_A(T))^{op}$.
\end{lem}
\noindent{\bf Proof.} Let $0\rightarrow T\rightarrow I_0\rightarrow
\cdots \rightarrow I_r\rightarrow 0$ be a minimal injective
resolution of $T$. Then we have the exact sequence $0\rightarrow
kQ\otimes_k T\rightarrow kQ\otimes_k I_0\rightarrow \cdots
\rightarrow kQ\otimes_k I_r\rightarrow 0$. By Lemma \ref{inj} ${\rm
inj.dim} (kQ\otimes_k I_j) \leq 1, \ \ 0 \leq j\leq r$, it follows
that ${\rm inj.dim} (kQ\otimes_k T)\leq r+1$.

Since the branch $(kQ\otimes_k T)_i$ is a direct sum of copies of
$T$, by Lemma \ref{adj3} we have
\begin{align*} {\rm Ext}_\m^s(kQ\otimes_k T, kQ\otimes_k T) & =
\bigoplus\limits_{i\in Q_0} {\rm Ext}_\m^s({\bf m}_i(T), kQ\otimes_k
T) \\ & \cong \bigoplus\limits_{i\in Q_0} {\rm Ext}_A^s(T,
(kQ\otimes_k T)_i) =0, \ \ \forall \ s\geq 1.\end{align*}

Now, put $\mathcal X = {\rm add} (kQ\otimes_k T)$. To see that
$kQ\otimes_k T$ is a cotilting $\m$-module, it remains to claim
$D(\m_\m)\in \widehat{\mathcal X}$, \ i.e., \ $D(kQ_{kQ})\otimes_k
D(A_A)\in \widehat{\mathcal X}$. In fact, since proj.dim
$D(kQ_{kQ})= 1$, we have an exact sequence $0\rightarrow
P_1\rightarrow P_0\rightarrow D(kQ_{kQ})\rightarrow 0$ with $P_0,
P_1$ being projective $kQ$-modules. So  we have the exact sequence
$0\rightarrow P_1\otimes_k D(A_A)\rightarrow P_0\otimes_k
D(A_A)\rightarrow D(kQ_{kQ})\otimes_k D(A_A)\rightarrow 0$. Since
$T$ is a cotilting $A$-module, we have an exact sequence
$0\rightarrow T_m\rightarrow \cdots \rightarrow T_0\rightarrow
D(A_A)\rightarrow 0$  with each $T_j\in {\rm add} (T)$. So we have
the exact sequence $0\rightarrow P_i\otimes_k T_m\rightarrow \cdots
\rightarrow P_i\otimes_kT_0\rightarrow P_i\otimes_kD(A_A)\rightarrow
0$, where $i = 0, 1$,  with each $P_i\otimes T_j\in {\rm add}
(kQ\otimes_k T)$. Thus $P_0\otimes_k D(A_A)\in \widehat{\mathcal X}$
and $P_1\otimes_k D(A_A)\in \widehat{\mathcal X}$. By Lemma \ref
{self-orthogonal}${\rm (i)}$ we have $D(kQ_{kQ})\otimes_k D(A_A)\in
\widehat{\mathcal X}$.

Finally, by Lemma \ref{adj2} we have
$${\rm Hom}_\m({\bf m}_i(T), {\bf m}_j(T))\cong {\rm Hom}_A(T, ({\bf m}_j(T))_i) = (\End_A(T))^{m_{ji}},$$
where  $m_{ji}$ is the number of paths of $Q$ from $j$ to $i$. Thus
one can easily see that there is an algebra isomorphism
$$\End_\m(kQ\otimes_k T)\cong \bigoplus\limits_{i, j\in Q_0}{\rm Hom}_\m({\bf m}_i(T), {\bf
m}_j(T)) \cong  (kQ\otimes_k \End_A(T))^{op}.$$ (In fact, if we
label the vertices of $Q$ as $1, \cdots, n$, such that if there is
an arrow from $j$ to $i$ then $j>i$. Then

$$kQ \cong \left(\begin{smallmatrix}
k&k^{m_{21}}&k^{m_{31}}&\cdots&k^{m_{n1}}\\
0&k&k^{m_{32}}&\cdots&k^{m_{n2}}\\
0&0&k&\cdots &k^{m_{n3}}\\
\vdots&\vdots&\vdots& &\vdots\\
0&0&0&\cdots&k\\
\end{smallmatrix}\right)_{n\times n}, $$
and hence

$$kQ\otimes_k \End_A(T)\cong \left(\begin{smallmatrix}
\End_A(T)&\End_A(T)^{m_{21}}&\End_A(T)^{m_{31}}&\cdots&\End_A(T)^{m_{n1}}\\
0&\End_A(T)&\End_A(T)^{m_{32}}&\cdots&\End_A(T)^{m_{n2}}\\
0&0&\End_A(T)&\cdots &\End_A(T)^{m_{n3}}\\
\vdots&\vdots&\vdots& &\vdots\\
0&0&0&\cdots&\End_A(T)\\
\end{smallmatrix}\right)_{n\times n}.)$$
This completes the proof. $\s$

\vskip10pt

\subsection{\bf Proof of Theorem \ref {mainthm2}.}
By Lemma \ref{cotilting} $kQ\otimes_k T$ is a cotilting $\m$-module,
and by Theorem \ref{mainthm1} ${\rm Mon}(Q, \ ^\perp T)= \ ^\perp
(kQ \otimes_k T)$.

If $L$ is another cotilting $\m$-module such that $^\perp L = {\rm
Mon}(Q, \ ^\perp T) = \ ^\perp (kQ \otimes_k T)$, then
$$\Ext_\m^s((kQ\otimes_k T)\oplus L, (kQ\otimes_k T)\oplus L) = 0, \ \ \forall \ s\ge 1,$$
so $(kQ\otimes_k T)\oplus L$ is also a cotilting $\m$-module. By [H]
the number of pairwise non-isomorphic direct summands of
$(kQ\otimes_k T)\oplus L$ is equal to the one of $kQ\otimes_k T$,
from which the proof is completed. $\s$

\vskip10pt

\section {\bf Contravariantly finiteness of monomorphism
categories}

\subsection{} \ Let $\mathcal X$ be a full subcategory of
$A$-mod and $M\in A$-mod. Recall from [AR] that {\it a right
$\mathcal X$-approximation of $M$} is an $A$-map $f:
X\longrightarrow M$ with $X\in\mathcal X$, such that the induced
homomorphism ${\rm{Hom}}_A (X', X)\longrightarrow {\rm{Hom}}_A(X',
M)$ is surjective for $X'\in\mathcal X$. If every $A$-module $M$
admits a right $\mathcal X$-approximation, then $\mathcal X$ is {\it
contravariantly finite in} $A$-mod. Dually one has the concept of
{\it a covariantly finite subcategory}. If $\mathcal X$ is both
contravariantly and covariantly finite, then $\mathcal X$ is {\it
functorially finite in} $A$-mod. Due to H. Krause and $\O$. Solberg
[KS, Corollary 0.3], a resolving contravariantly finite subcategory
is functorially finite, and a coresolving covariantly finite
subcategory is functorially finite. Due to M. Auslander and S. O.
Smal\o \ [AS, Theorem 2.4], a functorially  finite subcategory which
is closed under taking extensions has Auslander-Reiten sequences.

\subsection{} \ Auslander-Reiten [AR, Theorem 5.5(a)] claim that $\mathcal X$ is resolving
and contravariantly finite with $\widehat{\mathcal X} = A$-mod if
and only if $\mathcal X = \ ^\perp T$ for some cotilting $A$-module
$T$, where $\widehat{\mathcal X}$ is defined in 4.3.

\vskip10pt

As an application of Theorem \ref {mainthm2} and [AR, Theorem
5.5(a)], we have

\vskip10pt

\begin{thm} \label{app1} \ Let $\mathcal X$ be a full subcategory of $A$-mod. Then $\rm{Mon}(Q,\mathcal X)$ is a resolving contravariantly finite
 subcategory in $\m$-mod with $\widehat{{\rm Mon}(Q, \mathcal X)}= \m$-mod if and only if $\mathcal X$ is a resolving contravariantly finite subcategory
 in $A$-mod with $\widehat {\mathcal X} = A$-mod.
 \vskip5pt

 In particular, ${\rm Mon}(Q, A)$ is functorially
finite in ${\rm Rep}(Q, A)$, and ${\rm Mon}(Q, A)$ has
Auslander-Reiten sequences.
\end{thm}
\noindent{\bf Proof.} \ If $\mathcal X$ is resolving and
contravariantly finite with  $\widehat{\mathcal X}=A$-mod, then by
[AR, Theorem 5.5(a)] there is a cotilting module $T$ such that
$\mathcal X= \  ^\perp T$. By Theorem \ref {mainthm2} $kQ\otimes_k
T$ is a cotilting $\m$-module and ${\rm Mon}(Q, \mathcal X) = {\rm
Mon}(Q, \ ^\perp T) = \ ^\perp (kQ\otimes_k T)$, again by [AR,
Theorem 5.5(a)] we know that ${\rm Mon}(Q, \mathcal X)$ is resolving
and contravariantly finite with $\widehat{{\rm Mon}(Q, \mathcal X)}
= \m$-mod.

 Conversely, assume that ${\rm Mon}(Q, \mathcal X)$ is resolving and contravariantly finite with
 $\widehat {{\rm Mon}(Q,\mathcal X)} = \m$-mod. By Lemma \ref {exact2} $\mathcal X$ is resolving. To see that $\mathcal X$ is contravariantly finite,
 we take a sink in $Q_0$, say vertex $1$, and consider functor ${\bf m}_1: A\mbox{-}{\rm mod}\rightarrow {\rm Mon}(Q,
 A)$ (cf. 4.2). For $M\in A$-mod, since $1$ is a sink, ${\bf
 m}_1(M)$ has only one non-zero branch and its $1$-st branch is just
 $M$. Let $f: X\longrightarrow {\bf
 m}_1(M)$ be a right ${\rm Mon}(Q, \mathcal X)$-approximation. Then $f_1: X_1
 \longrightarrow M$ is a right $\mathcal X$-approximation (one can easily see this, for example, by Lemma \ref{adj2}. We omit the details).
By the same argument we see $\widehat{\mathcal
 X}=A$-mod since $\widehat {{\rm Mon}(Q,\mathcal X)} = \m$-mod. This
 completes the proof. $\s$

\vskip10pt

\section {\bf Finiteness of monomorphism categories}

As an application of Theorem \ref {mainthm2} and Auslander's
classical idea [Au, Chapter III], we describe the monomorphism
categories which are of finite type.

\vskip10pt

\subsection{} An additive full subcategory $\mathcal X$
of $A$-mod, which is closed under direct summands, {\it is of finite
type} if there are only finitely many isomorphism class of
indecomposable $A$-modules in $\mathcal X$.

\vskip10pt

An $A$-module $M$ is {\it an $A$-generator} if each projective
$A$-module is in ${\rm add} (M)$. A $\m$-generator $M$ is {\it a
generator and relative cogenerator of ${\rm Mon}(Q, A)$} if $M\in
{\rm Mon}(Q, A)$ and $kQ\otimes_k D(A_A)\in {\rm add} (M)$.

\vskip10pt

\begin{thm} \label{app2} ${\rm Mon}(Q, A)$ if of finite type if and only if there is a generator and relative  cogenerator $M$ of ${\rm Mon}(Q, A)$
 such that ${\rm gl.dim} \End_\m (M)\leq 2$.
 \end{thm}

\subsection{} \  Let $M$ be a $\m$-module. For an
arbitrary $\m$-module  $X$, denote by $\Omega_M(X)$ the kernel of a
minimal right ${\rm add}(M)$-approximation $M^\prime\longrightarrow
X$ of $X$. Define $\Omega_M^0(X)=X$, and
$\Omega_M^i(X)=\Omega_M(\Omega_M^{i-1}(X))$ for $i\ge 1$. Define
${\rm rel.dim}_M X$ to be the minimal non-negative integer $d$ such
that $\Omega_M^d(X)\in {\rm add}(M)$, or $\infty$ if otherwise. The
following fact is well known.

\vskip10pt

\begin{lem} {\rm(M. Auslander)} \ Let $M$ be an
$A$-module with $\Gamma=({\rm End}_A(M))^{op}$. Then for each
$A$-module $X$ we have ${\rm proj.dim}_\Gamma{\rm Hom}_A(M, X)\le
{\rm rel.dim}_M X$. Furthermore, if $M$ is a generator, then
equality holds.
\end{lem}

\vskip10pt

For an $A$-module $T$, denote by ${\mathcal X}_T$ the full
subcategory of $A$-mod given by
$$\{X\ |\ \exists \ \mbox{an exact sequence}\  0\rightarrow X\rightarrow T_0\xrightarrow{d_0}T_1\xrightarrow{d_1}\cdots,\ \mbox{with}\ T_i\in {\rm add}(T), \ \Ker d_i\in\ ^\perp T, \forall \ i\ge 0\}.$$
\noindent Not that ${\mathcal X}_T\subseteq \ ^\perp T$, and
${\mathcal X}_T= \ ^\perp T$ if $T$ is a cotilting module ([AR,
Theorem 5.4(b)]).

\vskip10pt

\begin{lem} \  Let $M$ be an $A$-generator with
$\Gamma=({\rm End}_A(M))^{op}$, and $T\in {\rm add}(M)$. Then for
each $A$-module $X\in {\mathcal X}_T$ and $X\not\in {\rm add}(T)$,
there is a $\Gamma$-module $Y$ such that ${\rm proj.dim}_\Gamma
Y=2+{\rm proj.dim}_\Gamma{\rm Hom}_A(M,X)$.
\end{lem}
\noindent {\bf Proof.} \ This is well-known. For completeness we
include a proof. By $X\in \mathcal{X}_{T}$ there is an exact
sequence $0 \rightarrow X \stackrel{u}\rightarrow T_0
\stackrel{v}\rightarrow T_1$ with $T_0, \ T_1 \in {\rm
add}(T)\subseteq {\rm add}(M)$. This yields an exact sequence
$$0\longrightarrow \Hom_A(M, X) \stackrel{u_*}\longrightarrow \Hom_A(M,
T_0) \stackrel{v_*}\longrightarrow \Hom_A(M, T_1) \longrightarrow
\cok v_*\longrightarrow 0.$$ Note the image of $v_*$ is not
projective (otherwise, $u_*$ splits, then we have an $A$-map $u':
T_0 \longrightarrow X$ such that $\Hom_A(M, u'u) = \Hom_A(M, {\rm
Id}_X)$. Since $M$ is an $A$-generator, we then get $u'u = {\rm
Id}_X$. This contradicts with $X\in {\rm add}(T)$). Putting $Y=\cok
v_*$,  we have ${\rm proj.dim}_\Gamma Y = 2+ {\rm proj.dim}_\Gamma
\Hom_A(M, X).$ \hfill $\s$

 \vskip10pt

\subsection{\bf Proof of Theorem \ref {app2}.} \ This is same as [Z, 5.3]. For completeness we include a proof.

Assume that ${\rm Mon}(Q, A)$ is of finite type. Then there is a
$\m$-module $M$ such that ${\rm Mon}(Q, A)={\rm add}(M)$. Since $kQ
\otimes_k D(A_A)\in{\rm Mon}(Q, A)$, and ${\rm Mon}(Q, A)$ contains
all the projective $\m$-modules, by definition $M$ is a generator
and relative cogenerator of  ${\rm Mon}(Q, A)$. Put $\Gamma= ({\rm
End}_A(M))^{op}$. For every $\Gamma$-module $Y$, take a projective
presentation ${\Hom}_\m(M, M_1)\xrightarrow{f_*}{\Hom}_\m(M,
M_0)\rightarrow Y\rightarrow 0$ of $Y$, where $M_1, M_0\in {\rm
add}(M)$, and $f: M_1\rightarrow M_0$ is a $\m$-map. Since ${\rm
inj.dim} (kQ\otimes D(A_A)) =1$ (Lemma \ref{inj}) and $M_1\in {\rm
Mon}(Q, A)=\ ^\perp (kQ\otimes D(A_A))$ (Proposition \ref{perp2}),
it follows that $\Ker f \in \ ^\perp(kQ\otimes D(A_A))= {\rm add}(M)
$. Thus
$$0\rightarrow{\Hom}_\m(M, \Ker f)\rightarrow{\Hom}_\m(M, M_1)\rightarrow{\Hom}_\m(M, M_0)\rightarrow Y\rightarrow 0$$
is a projective resolution of $\Gamma$-module of $Y$, i.e., ${\rm
proj.dim}_\Gamma Y\le 2$. This proves ${\rm gl.dim}_\Gamma Y\le 2$,
and hence ${\rm gl.dim}\End_\m(M)={\rm gl.dim}\Gamma\le 2$.

Conversely, assume that there is a generator and relative
cogenerator $M$ of ${\rm Mon}(Q,A)$ such that ${\rm gl. dim
}\End_\m(M)\le 2$. Put $\Gamma= ({\rm End}_A(M))^{op}$. Then ${\rm
gl.dim}\Gamma\le 2$. We claim that $ {\rm add}(M)= \
^\perp(kQ\otimes D(A_A))$, and hence by Proposition \ref{perp2}
${\rm Mon}(Q,A) = \ ^\perp(kQ\otimes D(A_A)) =  {\rm add}(M)$, i.e.,
${\rm Mon}(Q, A)$ is of finite type. In fact, since $ M \in{\rm
Mon}(Q,A) = \ ^\perp(kQ\otimes D(A_A))$, it follows that $ {\rm
add}(M) \subseteq \ ^\perp(kQ\otimes D(A_A))$. On the other hand,
let $X \in \ ^\perp(kQ\otimes D(A_A))$. By Theorem \ref{mainthm2}
 $kQ\otimes D(A_A)$ is a cotilting $\m$-module, and hence $^\perp(kQ\otimes D(A_A)) = \mathcal{X}_{kQ\otimes D(A_A)}$, by
[AR, Theorem 5.4(b)]. We divide into two cases. If $ X\in {\rm
add}(kQ\otimes D(A_A))$, then $X\in {\rm add}(M)$ since by
assumption $kQ\otimes D(A_A)\in {\rm add}(M)$. If $X \not\in {\rm
add}(kQ\otimes D(A_A))$, then by Lemma 6.3 there is a
$\Gamma$-module $Y$ such that ${\rm proj.dim}_\Gamma Y = 2 + {\rm
proj.dim}_\Gamma\Hom_\m(M,X)$. Now by Lemma 6.2 we have
$${\rm rel.dim}_M X={\rm proj.dim}_\Gamma\Hom_\m(M,X)={\rm
proj.dim}_\Gamma Y-2 \le 0,$$ this means $X \in {\rm add}(M)$. This
proves the claim and hence completes the proof. $\s$


\begin{thebibliography}{99}

\bibitem[A]{A} D. M. Arnold, Abelian groups and representations of finite partially
ordered sets, CMS Books Math., Springer-Verlag, New York, 2000.

\bibitem[Au]{Au} M. Auslander, Representation dimension of
artin algebras, Queen Mary College Math. Notes, Queen Mary College,
London, 1971; also in: I. Reiten, S. Smal${\o}$, ${\O}$. Solberg
(Eds.), Selected works of Maurice Auslander, Part 1 (II), Amer.
Math. Soc., 1999, 505-574.

\bibitem[AB]{AB}  M. Auslander, M. Bridger, Stable module theory, Mem. Amer. Math. Soc. vol. 94., Amer. Math. Soc., Providence, R.I., 1969.

\bibitem[AR]{AR} M. Auslander, I. Reiten, Applications of
convariantly finite subcategories, Adv. Math. 86(1991) 111-152.

\bibitem[AS]{AS} M. Auslander, S. O. Smal${\o}$, Almost split
sequences in subcategories, J. Algebra 69(1981) 426-454.

\bibitem[B]{B} G.
Birkhoff, Subgroups of abelian groups, Proc. Lond. Math. Soc. II,
Ser. (2) 38 (1934) 385-401.

\bibitem[C1]{C1} X. W. Chen, The stable
monomorphism category of a Frobenius category, Math. Res. Lett.
18(1)(2011) 125-137.

\bibitem[C2]{C2} X. W. Chen, Three results on Frobenius categories, Math. Z. 270 (1-2) (2012) 43 - 58.


\bibitem[EJ]{EJ} E. E. Enochs, O. M. G. Jenda, Relative
homological algebra, de Gruyter Exp. Math. vol. 30. Walter De
Gruyter Co., 2000.

\bibitem[H]{H} D. Happel, Triangulated categories in
representation theory of finite dimensional algebras, London Math.
Soc. Lecture Notes Ser. vol. 119, Cambridge Uni. Press, 1988.

\bibitem[HR]{HR} D. Happel, C.M.Ringel, Tilted algebras, Trans. Amer. Math. Soc. 274 (2) (1982) 399-443.

\bibitem[KLM1]{KLM1} D.
Kussin, H. Lenzing, H. Meltzer, Nilpotent operators and weighted
projective lines, to appear in: J. Reine Angew. Math. (Doi No.
10.1515/crelle-2012-0014),  avaible in arXiv: 1002.3797.

\bibitem[KLM2]{KLM2} D.
Kussin, H. Lenzing, H. Meltzer, Triangle singularities, ADE-chains,
and weighted projective lines,  to appear in: Adv. Math.,  avaible
in arXiv: 1203.5505.

\bibitem[KS]{KS} H.
Krause, O. Solberg, Applications of cotorsion pairs, J. London Math.
Soc. (2) 68 (3) (2003) 631-650.

\bibitem[LZ]{LZ} X. H. Luo, P. Zhang, Monic representations and Gorenstein-projective modules,
availble in arXiv 1110.6021.

\bibitem[Mi]{M} T. Miyashita, Tilting modules of finite projective dimension, Math.
Z. 193 (1986) 113-146.

\bibitem[Mo]{Mo} A. Moore, The Auslander and Ringel-Tachikawa theorem for submodule
embeddings, Comm. Algebra 38(2010) 3805-3820.

\bibitem[RW]{RW} F.
Richman, E. A. Walker, Subgroups of $p^5$-bounded groups, in:
Abelian groups and modules, Trends Math., Birkh\"auser, Basel, 1999,
pp. 55-73.

\bibitem[RS1]{RS1} C. M. Ringel, M. Schmidmeier, Submodules
categories of wild representation type, J. Pure Appl. Algebra
205(2)(2006) 412-422.

\bibitem[RS2]{RS2}
C. M. Ringel, M. Schmidmeier, The Auslander-Reiten translation in
submodule categories, Trans. Amer. Math. Soc. 360(2)(2008) 691-716.

\bibitem[RS3]{RS3}
C. M. Ringel, M. Schmidmeier, Invariant subspaces of nilpotent
operators I, J. rein angew. Math. 614 (2008) 1-52.

\bibitem[RZ]{RZ}
C. M. Ringel, P. Zhang, Representations of quivers over the algebras
of dual numbers, availble in arXiv 1112.1924.

\bibitem[S]{S}
D. Simson, Representation types of the category of subprojective
representations of a finite poset over $K[t]/(t^m)$ and a solution
of a Birkhoff type problem, J. Algebra 311(2007) 1-30.

\bibitem [SW]{SW} D. Simson,  M. Wojewodzki, An algorithmic solution of a
Birkhoff type problem, Fundamenta Informaticae 83(2008)  389-410.

\bibitem[XZZ]{XZZ} B. L. Xiong, P. Zhang, Y. H. Zhang,
Auslander-Reiten translations in monomorphism categories, Forum
Math. (DOI No. 10.1515/forum-2011-0003), avaible in arXiv:
1101.4113.

\bibitem[Z]{Z} P. Zhang, Monomorphism categories, cotilting
theory, and Gorenstein-projective modules, J. Algebra 339(2011)
180-202.
\end{thebibliography}
\end{document}